\documentclass[leqno,12pt]{amsart}

\usepackage{amsmath,amssymb,amsthm,amsfonts,amscd}
\usepackage[dvips]{graphics}

\newtheorem{dummy}{anything}
\newtheorem{theorem}[dummy]{Theorem}
\newtheorem{lemma}[dummy]{Lemma}

\newtheorem{corollary}[dummy]{Corollary}

\newtheorem{example}[dummy]{Example}

\newcommand{\Z}{{\mathbf Z}}
\newcommand{\N}{{\mathbf N}}

\newcommand{\tc}{{\sf {TC}}}
\newcommand{\cat}{\sf {cat}}
\newcommand{\RP}{{\mathbf {RP}}}
\newcommand{\cd}{{\sf {cd}}}

\newcommand{\supp}{{\sf {{supp}}}}
\newcommand{\cl}{{\rm cl}}

\begin{document}

\date{January 6, 2011}

\title[Random graph groups]{Topology of Random\\ Right Angled Artin Groups}

\author[Armindo ~Costa and Michael ~Farber]{Armindo ~Costa and Michael ~Farber}
\address{Department of Mathematical Sciences, Durham University\\
South Road, Durham, DH1 3LE, UK}

\email{a.e.costa@durham.ac.uk}
\email{michael.farber@durham.ac.uk}

\maketitle
\begin{abstract} In this paper we study topological invariants of a class of random groups. Namely, we study right angled Artin groups associated to random graphs
and investigate their Betti numbers, cohomological dimension and topological complexity. The latter is a numerical homotopy invariant reflecting complexity
of motion planning algorithms in robotics. We show that the topological complexity of a random right angled Artin group assumes,
with probability tending to one,
at most three values,  when $n\to \infty$. We use a result of Cohen and Pruidze which expresses the topological complexity of right angled Artin groups in combinatorial terms. 
Our proof deals with the existence of bi-cliques in random graphs.
\end{abstract}

\section{Introduction}

Problems of mathematical modeling of large systems of various nature (economical, mechanical, ecological, etc)
motivate studying random geometric, topological and algebraic objects. 
For a system of great complexity, it is unrealistic to assume that one may have a precise description of its configuration space; it is more reasonable to view it as a 
partially known or a random space.

Studying random topological and algebraic objects instead of their deterministic
analogues has several major advantages. Firstly, random mathematical objects in many cases more adequately reflect reality. Secondly, random objects are often simpler mathematically since exotic and most complicated examples can be ignored if they are rare, i.e. appear with very small probability. 

Finally, the probabilistic method is useful in proving existence theorems and in building new objects and counter-examples. 
Recently Gromov \cite{G} used the theory of random groups to prove the existence of a finitely presentable group whose Cayley graph 
contains a family of expanding graphs.

The most developed stochastic--topological object is a random graph. 
The theory of random graphs, initiated circa 1959 by Erd\"os and R\'enyi \cite{ER}, is nowadays a fast growing branch of applied mathematics, 
offering a plethora of spectacular results and predictions for various engineering and computer science applications, see \cite{AS}, \cite{B}, \cite{JLR}. 

Random simplicial complexes of high dimension were recently suggested and studied by Linial--Meshulam in~\cite{LM}, and 
Meshulam--Wallach in~\cite{MW}. 
The fundamental groups of random 2-complexes are random groups of a fairly general type: clearly this model covers all finitely presented groups. 
In a recent paper E. Babson, C. Hoffman and M. Kahle \cite{BHK} show that a random group in this model is hyperbolic (or trivial)\footnote{The result of \cite{BHK} leaves open the possibility that a random group in the Linial-Meshulam model with probability parameter $p$ near $n^{-1/2}$ is not hyperbolic.}. 
This result is of great interest since it is known that many 
delicate  problems of group theory are algorithmically solvable for hyperbolic groups. 
Paper \cite{CFK} establishes a threshold for freeness of a random group produced by the Linial - Meshulam model. 

Configuration spaces of mechanical linkages with random bar lengths were studied in papers
\cite{F1}, \cite{FK}. These are closed smooth manifolds depending on a large number of independent random parameters. Although the number of homeomorphism types of these manifolds grows extremely fast with dimension, their topological characteristics can be predicted with high probability when the number of links tends to infinity.

The theory of random groups was initiated by M. Gromov \cite{Gr}, \cite{G} who proposed several different models generating random groups. 
One of this models (called the density model) has a density parameter $0\le d\le 1$. A random group in this model is given by a presentation
$G=\langle a_1, \dots, a_m|\, R\, \rangle$
where the number of generators $m\ge 2$ is fixed and the set of relations $R$ consists of $(2m-1)^{d\ell}$ words which are picked randomly from the set 
$S_\ell$ of all reduced words in $a_1, \dots, a_m$ of length $\ell$. 
One is interested in probabilities of various properties of these groups as the length parameter $\ell$ tends to infinity. One says that a group property holds with overwhelming probability
if the probability of its occurence tends to 1 as $\ell \to \infty$. 
The striking phase transition theorem of Gromov states 
that in the density model, with overwhelming probability,
for $d<1/2$ 
a random group is infinite, hyperbolic, torsion free, of geometric dimension 2, 
while for $d>1/2$ it is either $\Z_2$ or trivial; see \cite{Gr}, pp. 273-275 and \cite{Oll}, page 31. 
We refer also to  \cite{Oll} for more details. 

In \cite{Zuk} A. Zuk studied an interesting model of random groups (the triangular model) where all relations have length 3 but the number of generators $m$ is not fixed and tends to infinity. 
We refer to \cite{Oll}, page 41 where this model is compared with the Gromov density model. 

The article \cite{S} is a recent survey of probabilistic group theory.  We also refer to a recent paper \cite{KS} which compares different models of random groups.

Given a finite graph (i.e. a one-dimensional simplicial complex) $\Gamma$ with vertex set $V$ and with the set of edges $E$ one associates to it a right angled Artin group (RAAG) (also known as a graph group)
$$G_\Gamma= \langle v\in V; vw=wv \quad \mbox{iff} \quad (v,w)\in E\rangle,$$
see \cite{C}, \cite{MV}.
In the case when $\Gamma$ is a complete graph $G_\Gamma$ is a free abelian group of rank $n=|V|$; in the other extreme, when $\Gamma$ has no edges the group $G_\Gamma$ is the free group of rank $n$. In general $G_\Gamma$ interpolates between the  free and free abelian groups.
In this paper we are interested in right angled Artin groups associated to {\it random graphs} $\Gamma$. We adopt one of the basic Erd{\H{o}}s
- R\'{e}nyi models of random graphs in which each edge of the complete graph on $n$ vertices is included with probability $0<p<1$ independently of all other edges. In other words, we consider the probability space $\Omega_n$ of all $2^{\binom n 2}$ subgraphs of the complete graph on $n$ vertices $\{1, 2, \dots,n\}$ and the probability that a specific graph $\Gamma\in \Omega_n$ appears as a result of a random process equals
\begin{eqnarray}\label{prob}
{\mathbf P}(\Gamma)\,  = \, p^{E_\Gamma}(1-p)^{{\binom n 2}-E_\Gamma},
\end{eqnarray}
where $E_\Gamma$ denotes the number of edges of $\Gamma$, see \cite{JLR}.

In Gromov's model a random group is given by a presentation with randomly generated relations where each of the relations has a fixed length $\ell$
and the total number of relations is also fixed (and depends on the density $d$), see above. In the case of random right angled Artin groups one generates  randomly relations of a special type (namely, the commutation relations).  


In this paper we examine statistics of various topological invariants of the group $G_\Gamma$ associated to a random graph.
Each of such invariants is a random function and it is quite natural to ask about its mathematical expectation and distribution function. We assume that $n\to \infty$ and seek results of asymptotic nature. 

Our main effort is in finding the topological complexity $\tc(K_\Gamma)$ of the Eilen\-berg-MacLane spaces $K_\Gamma$ associated to the right angled Artin group of a random graph $\Gamma$. 
The spaces $K_\Gamma$ are closely related to the configuration spaces of braid groups of graphs (see \cite{CW}) which indicates the relevance of the results obtained here to collision free motion planning of multiple particles on graphs.  The notion of topological complexity $\tc(X)$ reflects the structure of motion planning algorithms for systems having $X$ as their 
configuration space, see a brief survey and further references in \S 4. In general, $\tc(X)$ depends only on the homotopy type of $X$; however if $X=K(G,1)$ is an Eilenberg - MacLane space then $\tc(X)$ depends only on the fundamental group $G$ of $X$. In this paper we show that the topological complexity $\tc(K_\Gamma)$, where $\Gamma$ is a random graph, can be determined almost precisely (with ambiguity at most two). We prove the following result:
\vskip 0.5cm
{\bf Theorem. }{\it Fix an arbitrary $\epsilon >0$ and assume that the probability parameter $0<p<1$ is constant (i.e. independent of $n$). 
Then for a random graph $\Gamma\in \Omega_n$ one has
\begin{eqnarray}\label{ineq122}
2\cdot \lfloor z(n,p)-\epsilon\rfloor +1 \, \leq\,  \tc(K_\Gamma)\, \leq \, 2\cdot \lfloor z(n,p) +\epsilon\rfloor+1,
\end{eqnarray}
asymptotically almost surely, where 
$$z(n, p) =   2\log_q n - 2\log_q \log_q n + 2 \log_q(e/2) +1, \quad q= p^{-1}. $$
In other words, probability that a graph $\Gamma\in \Omega_n$ does not satisfy inequality (\ref{ineq122}) tends to zero when $n$ tends to infinity.}

A recent paper \cite{ChF} continues investigation of properties of random groups associated to random graphs initiated by the present article. It is shown in \cite{ChF} that a random right angled Artin group has a finite outer 
automorphism group if the probability parameter $p$ is constant and lies in the interval $(0.3, 1)$. Paper \cite{ChF} also deals with a more general class of random groups associated to random graphs and gives a threshold for hyperbolicity of these groups. 

The authors thank Daniel Cohen for helpful discussions.

\section{Betti numbers of random graph groups}\label{betti}

First we remind the well-known construction of an aspherical complex $K_\Gamma$ with fundamental group $G_\Gamma$. We refer to \cite{C} and \cite{MV} for proofs and more detail.

Let $V=V_\Gamma$ denote the set of vertices of the graph
$\Gamma$. The torus $T^n$ where $n=|V|$ can be identified with the set of all functions $\phi: V \to S^1$.
The support $\supp(\phi)\subset V$ of a function $\phi: V \to S^1$  is defined as the set of vertices $v\in V$ such that $\phi(v) \not=1$.
One defines $K_\Gamma\subset T^n$ to be the set of all functions $\phi$ such that their support $\supp(\phi)$ generates a complete subgraph of $\Gamma$, i.e. any two vertices of the support
are connected by an edge in $\Gamma$. It is known \cite{C}, \cite{MV} that $K_\Gamma$ (viewed with the induced topology) is aspherical and its fundamental group is $\Gamma$.

Fix the cell decomposition of $S^1$ consisting of a 0-cell $1\in S^1$ and a 1-cell $S^1-\{1\}$. Then $T^n$ inherits a cell decomposition with cells in one-to-one correspondence with
subsets of $V$. In this decomposition $K_\Gamma \subset T^n$ is a cell subcomplex; the cells of $K_\Gamma$ are in 1-1 correspondence with complete subgraphs of $\Gamma$.
Namely, given a subset $S\subset V$ one considers the set $e_S$ of all functions $\phi: V \to S^1$ with support equal $S$; clearly $e_S$ is a cell of dimension $|S|$.

The cohomology algebra of $K_\Gamma$ with integral coefficients is the quotient
\begin{eqnarray}
H^\ast(K_\Gamma; \Z) \simeq E(v_1, \dots, v_n)/J_\Gamma\end{eqnarray}
where $E(v_1, \dots, v_n)$ is the exterior algebra generated by degree one classes corresponding to the vertices $V=\{v_1, \dots, v_n\}$ of $\Gamma$ and the ideal $J_\Gamma$ is generated by the degree two monomials
$vw$
such that the corresponding vertices $v, w$ are not connected by an edge.

In particular any product $v_{i_1} v_{i_2} \dots v_{i_r}$ vanishes iff the
corresponding vertices $\{v_{i_1}, v_{i_2}, \dots, v_{i_r}\}$ do not form a complete subgraph of $\Gamma$. The monomials of the form 
$v_{i_1} v_{i_2} \dots v_{i_r}$, with distinct vertices $v_{i_1}, v_{i_2}, \dots, v_{i_r}$ such that any two of them are connected by an edge in $\Gamma$,
form an additive basis of $H^\ast(K_\Gamma;\Z)$. Thus, one obtains the following well-known fact: 

\begin{lemma}\label{lmbetti}
For an integer $r\ge 2$ the $r$-th Betti number\footnote{The discussion preceeding Lemma \ref{lmbetti} shows that the cohomology $H^\ast(K_\Gamma;\Z)$ is $\Z$-torsion free. Thus the Betti numbers 
$b_r(G_\Gamma)$ are independent of the field of coefficients.}
 $b_r(G_\Gamma)=b_r(K_\Gamma)$ equals the number of complete subgraphs of size $r$ in $\Gamma$.
\end{lemma}
Note that $b_0(G_\Gamma)=1$ and $b_1(G_\Gamma)=n$ for any graph $\Gamma$ on $n$ vertices. 

\begin{lemma}\label{lmtwo}
The expectation of the $r$-th Betti number of the group $G_\Gamma$ of a random graph $\Gamma$, where $r\ge 2$, equals
\begin{eqnarray}\label{expectation}
{\rm \mathbf E}(b_r(G_\Gamma))={\binom n r}p^{\binom r 2}.
\end{eqnarray}
\end{lemma}
\begin{proof}
As explained above we must find the number of maximal complete subgraphs of size $r$ in $\Gamma\in \Omega_n$. For a subset $S\subset \{1, \dots, n\}$ with $|S|=r$ consider the random variable
$I_S : \Omega_n\to \{0,1\}$ which equals 1 on a graph $\Gamma\in \Omega_n$ iff $S$ forms a complete subgraph in $\Gamma$. Then ${\mathbf E}(I_S)= p^{\binom r 2}$ and
$\sum_S I_S$ is the number of all complete subgraphs on $r$ vertices. This shows that ${\mathbf E} (\sum_S I_S)$ is as stated.
\end{proof}

Now we assume that $r$ (the dimension) is fixed and $p$ may depend on $n$.
Asymptotically, the expectation of $b_r(G_\Gamma)$ can be written as

\begin{eqnarray*}{\rm \mathbf E}(b_r(G_\Gamma)) \sim \frac{1}{r!} \left[ np^{\frac{r-1}{2}}\right]^r.
\end{eqnarray*}
The expectation has a positive limit for $n\to \infty$ if and only if
\begin{eqnarray}\label{limit}
np^{\frac{r-1}{2}} \to c>0.
\end{eqnarray}
Under this condition the expectation ${\rm \mathbf E}(b_r(G_\Gamma))$ converges to $\displaystyle \frac{c^r}{r!}.$

Note that the convergence (\ref{limit}) to a positive limit may happen for one dimension $r$ only.

Moreover, assuming (\ref{limit}), the distribution of $b_r: \Omega \to \Z$ converges to the Poisson distribution with expectation
\begin{eqnarray}\label{lambda1}
\lambda = \displaystyle \frac{c^r}{r!},
\end{eqnarray}
see below. Theorem \ref{inter} is a group theoretic interpretation (based on Lemma \ref{lmbetti}) 
of a theorem of Sch\"urger \cite{Sch} about complete subgraphs in random graphs.
More general results concerning the containment of a specific graph in a random graph were established by by Bollob\'{a}s \cite{B1} and Karo\'{n}ski and Ruci\'{n}ski \cite{KR};
see
 Theorem 3.19 from \cite{JLR}.

\begin{theorem}\label{inter}  Fix an integer $r>1$ and consider the function of $r$-th Betti number of the associated graph group,
$$b_r: \Omega_n \to \Z, \quad b_r(\Gamma) = b_r(G_\Gamma),$$
as a random function of a random graph. If the limit (\ref{limit}) exists and is positive then for any integer $k=0, 1, \dots$ the probability
$${\mathbf P}(b_r(G_\Gamma)=k) $$
converges (as $n\to \infty$) to
$$e^{-\lambda}\cdot \displaystyle \frac{\lambda^k}{k!}$$ where $\lambda$ is given by formula (\ref{lambda1}).
\end{theorem}

In other words, Theorem \ref{inter} claims that the limiting distribution is Poisson with mean $\lambda$.

\begin{example} {\rm Consider the following examples illustrating the previous Theorem.

Suppose that $r=2$ and $p=\frac{4}{n^2}$. Then $c=2$, $\lambda=2$, and for any integer $k=0, 1, \dots$
the probability that $b_2(G_\Gamma)=k$ converges to $\frac{2^k}{e^2\cdot k!}$ as $n\to \infty$.

As another example, assume that $r=3$ and $p=\frac{6}{n}$. Then $\lambda =36$ and
the probability that $b_3(G_\Gamma)=k$ converges to $\frac{36^k}{e^{36}\cdot k!}$ as $n\to \infty$
}
\end{example}
\section{Cohomological dimension of random graph groups}

The cohomological dimension of $G_\Gamma$
equals the size of the maximal clique in $\Gamma$; this follows from the discussion of section \ref{betti}.
Recall that a clique in a graph is defined as a maximal complete subgraph. The clique number $\cl(\Gamma)$ of a graph $\Gamma$ is the maximal order
of a clique in $\Gamma$.

There are many results in the literature about the clique number of random graphs; we may interpret these results as statements about the cohomological dimension
of graph groups build out of random graphs.  Matula \cite{M1}, \cite{M2} discovered that for fixed values of $p$ the distribution of the clique number of a random graph is highly
concentrated in the sense that almost all random graphs have about the same clique number. These results were developed further by Bollob\'{a}s and Erd\H{o}s \cite{BE}; see the monographs of B. Bollob\'{a}s \cite{B} and of N. Alon and J. Spencer \cite{AS}.

Below we restate a result of Matula \cite{M2} as a statement about cohomological dimension of random graph groups.

Denote
\begin{eqnarray}\label{approx}
z(n, p) =   2\log_q n - 2\log_q \log_q n + 2 \log_q(e/2) +1,
\end{eqnarray}
where $q=p^{-1}$. We assume that $p$ is independent of $n$, i.e. it is constant.  

\begin{theorem}\label{matula} Fix an arbitrary $\epsilon >0$. Then
\begin{eqnarray}\label{ineq}
\lfloor z(n,p)-\epsilon\rfloor \, \leq \, \cd(G_\Gamma)\, \leq \, \lfloor z(n,p) +\epsilon\rfloor,
\end{eqnarray}
asymptotically almost surely (a.a.s). In other words, the probability that a graph $\Gamma\in \Omega_n$ does not satisfy inequality (\ref{ineq}) tends to zero when $n$ tends to infinity.
\end{theorem}

Here $\lfloor x \rfloor $ denotes the largest integer not exceeding $x$. We may assume that $\epsilon < 1/2$; then the integers
$\lfloor z(n,p)-\epsilon\rfloor$ and $\lfloor z(n,p) +\epsilon\rfloor$ either coincide or differ by $1$.

Thus, according Theorem \ref{matula}, the cohomological dimension $\cd(G_\Gamma)$ for a random graph $\Gamma$ takes on one of at most two values depending on $n$ and $p$, with probability approaching $1$ as $n \to \infty$.

Next Lemma is a technical result which will be used later in this paper. The Lemma shows in particular that the expected Betti number of $K_\Gamma$ 
in dimension $r=\lfloor z(n,p)-\epsilon\rfloor$ (nearly the highest dimension where homology is nonzero, according to Lemma \ref{lmtwo} and Theorem \ref{matula}) tends to infinity faster than $r$. Thus, significant amount of homology of $K_\Gamma$ is concentrated in the top dimension. 

\begin{lemma}\label{lemma4}
Fix $\epsilon >0$ and let  $r=\lfloor z(n,p)-\epsilon\rfloor$. Then
$$r^{-1}\cdot{ \binom n r}p^{\binom r 2}\to \infty$$
as $n\to \infty$.
\end{lemma}
\begin{proof} One has $r=\lfloor z(n,p) -\epsilon\rfloor \le z(n,p)-\epsilon$ and therefore
\begin{eqnarray*}
p^{\binom r 2} \ge \left( p^{(z(n,p) - \epsilon -1)/2}\right)^r
\\ =\left(p^{\log_q n- \log_q\log_q n +\log_q(e/2) - \epsilon/2}\right)^r\\
=\left(\frac{2C\log_q n}{en}\right)^r, \quad \mbox{where}\quad C=q^{\epsilon/2}>1.
\end{eqnarray*}
On the other hand, using Stirling's formula, we have
\begin{eqnarray*}
{\binom n r} = c_n \cdot\left(\frac{n}{r}\right)^r\cdot e^r \cdot r^{-1/2}
\end{eqnarray*}
where $c_n$ and $c_n^{-1}$ are bounded. Therefore,
\begin{eqnarray*}
r^{-1}\cdot {\binom n r}\cdot p^{\binom r 2} \ge r^{-1}c_n\left(\frac{n}{r}\right)^rr^{-1/2}\left(\frac{2C\log_qn}{en}\right)^r e^r=\\
\left(C\cdot \frac{2\log_q n}{r}\right)^r \cdot r^{-3/2}\cdot c_n \ge C^r\cdot r^{-3/2}\cdot c_n.
\end{eqnarray*}
Clearly, $C^r\cdot r^{-3/2}\cdot c_n$ tends to infinity since $C>1$. This completes the proof.
\end{proof}

\section{Motion planning algorithms and the concept of topological complexity}

Given a mechanical system, a motion planning algorithm is a function which assigns to any pair of states of the system, an initial state and a desired state, a continuous motion of the system starting at the initial state and ending at the desired state. The design of effective motion planning algorithms is one of
the challenges of modern robotics, see \cite{La}.
Motion planning algorithms are applicable in various situations when the system is autonomous and operates in a fully or partially known environment.

The complexity of motion planning algorithms is measured by a numerical invariant $\tc(X)$ which depends on the homotopy type of the configuration space
$X$ of the system \cite{F1}. This invariant is defined as the Schwarz genus (also known as the \lq\lq sectional category\rq\rq) of the path-space fibration
\begin{eqnarray}\label{fibration}
p: PX \to X\times X.\end{eqnarray}
Here $PX$ is the space of all continuous paths $\gamma: [0,1]\to X$ equipped with the compact-open topology and $p(\gamma)=(\gamma(0), \gamma(1))$ is the map
associating to a path $\gamma:[0,1]\to X$ its pair of endpoints. Explicitly, $\tc(X)$ is the smallest integer $k$ such that $X\times X$ admits an open cover
$U_1\cup U_2\cup \dots\cup U_k=X\times X$ with the property that there exists a continuous section $U_i\to PX$ of (\ref{fibration})  for each $i=1, \dots, k$.
If $X$ is an Euclidean neighborhood retract then $\tc(X)$ can be equivalently characterized as the minimal integer $k$ such that there exists a section $s: X\times X\to PX$ of the fibration $p$ with the property that $X\times X$ can be represented as the union of $k$ mutually disjoint locally compact sets
$$X\times X=G_1\cup \dots\cup G_k$$
such that the restriction $s|G_i$ is continuous for $i=1, \dots, k$, see \cite{invit}, Proposition 4.2. A section $s$ as above represents a motion planning algorithm: given a pair $(A,B)\in X\times X$ the image $s(A, B)\in PX$ is a continuous motion of the system starting at the state $A$ and ending at the state $B$.

Intuitively, the topological complexity $\tc(X)$ can be understood as a measure of the navigational complexity of the topological space $X$; it is the minimal number of continuous rules which are needed to describe a motion planning algorithm in $X$.

The invariant $\tc(X)$ admits an upper bound in terms of the dimension of the configuration space $X$,
\begin{eqnarray}\label{upperbound}
\tc(X) \le 2\dim(X) +1
\end{eqnarray}
see \cite{F1}, Theorem 4. There are many examples when inequality (\ref{upperbound}) is sharp: take for instance $X= T^n\sharp T^n$, the connected sum of two copies of a torus,
having the topological complexity $\tc(X)= 2n+1$. However for any simply connected space $X$ one has a more powerful upper bound
\begin{eqnarray}\label{ineq2}\tc(X) \le \dim(X) +1,\end{eqnarray}
see \cite{F2}. The latter inequality is sharp for any simply connected closed symplectic manifold $X$, see \cite{FTY}.

There are many other examples when
the inequality (\ref{upperbound}) is not sharp. It was established in \cite{CF} that for any finite cell complex $X$ with $\pi_1(X)=\Z_2$ one has
\begin{eqnarray}\label{upperbound1}\tc(X) \le 2\dim (X).\end{eqnarray}

For example $\tc(\RP^n)\le 2n$ for all $n$; moreover, $\tc(\RP^n)=2n$ if and only if $n$ is a power of $2$, see Corollary 14 of \cite{FTY}.

The main result of this paper states that the inequality (\ref{upperbound}) is asymptotically very close to equality in the case of Eilenberg - MacLane spaces of random graph groups.

\section{The topological complexity of random graph groups}

Consider the probability space $\Omega_n$ of random graphs on $n$ vertices with probability given by formula (\ref{prob}). For any $\Gamma\in \Omega_n$ consider the corresponding Eilenberg-MacLane complex $K_\Gamma$ (see \S 2) and its topological complexity $\tc(K_\Gamma)$, as defined in the previous section. 
The probability measure on $\Omega_n$ is given by the formula (\ref{prob}). 

\begin{theorem}\label{main1}
Fix an arbitrary $\epsilon>0$ and assume that the edge probability parameter $0<p<1$ is constant (i.e. independent of $n$). Then for any random graph $\Gamma\in \Omega_n$ one has
\begin{eqnarray}\label{ineq1}
2\cdot \lfloor z(n,p)-\epsilon\rfloor +1 \, \leq\,  \tc(K_\Gamma)\, \leq \, 2\cdot \lfloor z(n,p) +\epsilon\rfloor+1,
\end{eqnarray}
asymptotically almost surely. Here $z(n,p)$ is given by formula (\ref{approx}). 
\end{theorem}

It is clear that the integers on the left and on the right of inequality (\ref{ineq1}) differ at most by $2$ (if $\epsilon <1/2$). Hence Theorem \ref{main1} determines the value of the topological complexity $\tc(G_\Gamma)$ for a random graph with ambiguity of at most 2. Comparing with the result of Theorem \ref{matula} we obtain
\begin{corollary}
For a random graph $\Gamma\in \Omega_n$ one has
\begin{eqnarray}\label{equal}
 2\cdot \cd(G_\Gamma) -1 \le \tc(K_\Gamma)  \le 2\cdot \cd(G_\Gamma) +1,
\end{eqnarray}
asymptotically almost surely.
\end{corollary}

Note that we have $$\cd(G_\Gamma) =\dim K_\Gamma= \cat(K_\Gamma)-1$$ for the Lusternik - Schnirelmann category, as it is easy to see.

The rest is this section is devoted to the proof of Theorem \ref{main1}.

By an $(r,r)$ {\it bi-clique} in a graph $\Gamma$ we understand an ordered pair consisting of two vertex disjoint complete subgraphs of $\Gamma$, each having $r$ vertices.
To specify an $(r,r)$ {\it bi-clique} one has to determine an $r$-element subset $S$ of the set of vertices of $\Gamma$ and an $r$-element subset $T$ in the complement $V-S$ such that the induced graphs on $S$ and $T$ are complete.

We know from sections \S\S 2, 3 that $\cd(G_\Gamma) \ge r$ if $\Gamma$ contains an $r$-clique, i.e. a maximal complete subgraph on $r$ vertices.
By Theorem of Cohen and Pruidze \cite{CP} one has $\tc(K_\Gamma) \ge 2r+1$ if $\Gamma$ contains an $(r,r)$ bi-clique.

In the rest of this section we set
$$r= \lfloor z(n,p)-\epsilon\rfloor .$$ Theorem \ref{main1} follows once we have shown that a random graph $\Gamma\in \Omega_n$ contains an $(r,r)$ bi-clique a.a.s. The right hand side of the inequality (\ref{ineq1}) follows from the general upper bound (\ref{ineq2}) and from the right hand side of (\ref{ineq}).

Let $r>0$ be an integer and let $X: \Omega_n\to \Z$ denote the number of $(r,r)$ bi-cliques in random graph. Our goal is to show that $X>0$ asymptotically almost surely,
i.e.
\begin{eqnarray}\label{pos}
{\mathbf P}(X>0) \to 1, \quad \mbox{for}\quad n\to \infty.\end{eqnarray}
The proof of (\ref{pos}) will use the second moment method and will be based on the inequality
\begin{eqnarray}
{\mathbf P}(X>0) \, \ge \, \frac{({\mathbf E}X)^2}{{\mathbf E}(X^2)}
\end{eqnarray}
see \cite{JLR}, page 54. Thus, our statement follows once we show that
\begin{eqnarray}
\frac{{\mathbf E}(X^2)}{({\mathbf E}X)^2} \, \to \, 1\quad \mbox{as} \quad n\to \infty.
\end{eqnarray}

Let $S$ and $T$ be disjoint $r$-element subsets of the set of vertices of $K_n$ and let $$I_{(S,T)}: \Omega_n \to \{0,1\}$$ denote the function which equals 1 on a graph $\Gamma \in \Omega_n$ if and only if $S$ and $T$ form a bi-clique in $\Gamma$. Then
$$X=\sum_{(S, T)} I_{(S,T)}$$
where the sum is taken over all ordered pairs of disjoint $r$-element subsets of $\{1, 2, \dots, n\}$. Note that one obviously has
$${\mathbf E}(I_{(S, T)})\,  =\,  p^{2{\binom r 2}}$$
and thus
$${\mathbf E}(X) = {\binom n {r, r}} p^{2{\binom r 2}}\, \, ,$$
where
$$\binom n {r, r} = \frac{n!}{r! \cdot r! \cdot (n-2r)!}$$
denotes the multynomial coefficient. Similarly,
\begin{eqnarray}\label{square}
X^2 = \sum I_{(S,T)}\cdot I_{(S', T')} .\end{eqnarray} 
Here $(S, T)$ and $(S', T')$ run over all ordered pairs of disjoint $r$-element subsets of the set of vertices $\{1, \dots, n\}$.
\begin{figure}[h]
\begin{center}
\resizebox{7cm}{5cm}
{\includegraphics[100,400][476,724]{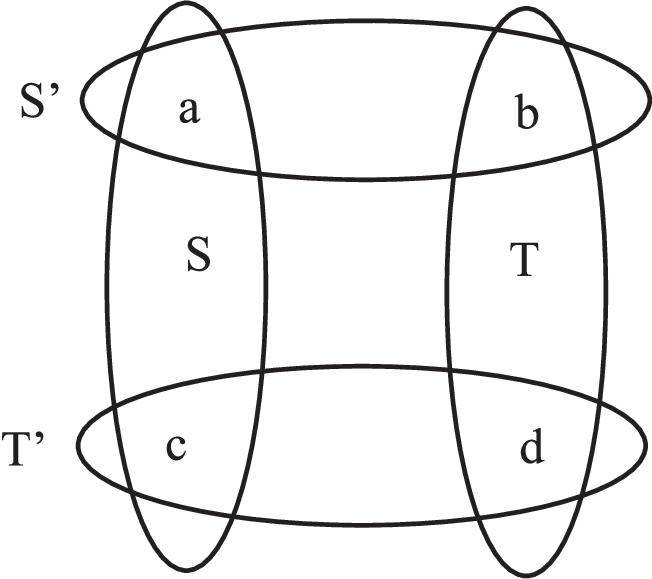}}
\end{center}
\caption{A pair of $(r,r)$ bi-cliques}\label{fig1}
\end{figure}
Denoting
$$a = |S\cap S'|, \, b= |T\cap S'|, \, c = |S\cap T'|, \, d = |T\cap T'|,$$
(see Figure \ref{fig1}) we find
\begin{eqnarray}
{\mathbf E}(I_{(S,T)}\cdot I_{(S', T')}) = {\displaystyle p^{4{\binom r 2} - {\binom a 2} -{\binom b 2}- {\binom c 2}- {\binom d 2}}}.
\end{eqnarray}
Therefore taking into account (\ref{square}) one obtains the following expression
\begin{eqnarray}\label{formula}
\frac{{\mathbf E}(X^2)}{{\mathbf E}(X)^2 }\, =\,  \sum_{\alpha\in D} \, F_\alpha \cdot q^{L(\alpha)}  = \, \sum_{\alpha\in D} \, T_\alpha.
\end{eqnarray}
Here $$\alpha=(a, b, c, d)\in \Z^4$$ denotes a vector and $D$ is the set of all vectors $\alpha=(a, b, c, d)$
with nonnegative integer components satisfying the inequalities
\begin{eqnarray}\label{domain}
a+b \le r, \quad a+c \le r, \quad c+d \le r, \quad b+d\le r,
\end{eqnarray}
see Figure \ref{fig1}.  
In formula (\ref{formula}) the coefficient $F_\alpha$
is given by
\begin{eqnarray}F_\alpha = {\displaystyle \frac{{\binom r {a,\, c}}{\binom r {b,\,  d}} {\binom {n-2r} {r-a-b, \, r- c-d}}}{\binom n {r, r}}}\end{eqnarray}
  and
\begin{eqnarray}L(\alpha) = {\binom a 2}+{\binom b 2}+{\binom c 2}+{\binom d 2},\quad \quad q =p^{-1},\end{eqnarray}
while \begin{eqnarray}T_\alpha = F_\alpha \cdot q^{L(\alpha)}.\end{eqnarray}

Let ${\mathbf m}(x,y)$ denotes $\max\{x, y\}$. Then inequalities (\ref{domain}) can be rewritten in a simple form as
\begin{eqnarray}
{\mathbf m}(a,d) + {\mathbf m}(b,c) \le r.
\end{eqnarray}

Next we mention the symmetry of the problem. There are two commuting involutions
\begin{eqnarray*}
\beta, \gamma: D\to D, \quad \beta^2=1=\gamma^2,
\end{eqnarray*}
where
$$\beta(a) =b, \quad \beta(c)=d, \quad \gamma(a) = c, \quad \gamma(b)=d.$$
These two involutions generate an action of the group $G=\Z_2\oplus \Z_2$ of $D$ which preserves both functions $T_\alpha$ and $L(\alpha)$.
This action is transitive on the four coordinates.

Recall that our goal is to show that the sum (\ref{formula}) tends to 1 as $n\to \infty$.
Note that
 \begin{eqnarray}
 \sum_{\alpha\in D} F_\alpha \, = \, 1
 \end{eqnarray}
for obvious reasons. Observe also that the term $F_0$ corresponding to $\alpha = (0,0,0,0)\in D$ equals
\begin{eqnarray*}
F_0 &=& \frac{\binom {n-2r} {r, \, r}}{\binom n {r, \, r}} = \prod_{k=0}^{2r-1} \left( 1- \frac{2r}{n-k}\right) \\
 &\geq& \left(1- \frac{2r}{n-2r+1}\right)^{2r}\ge 1 - \frac{4r^2}{n-2r+1}.\end{eqnarray*}
Hence we see that $F_0\to 1$ as $n\to \infty$. Therefore, the sum of all coefficients $F_\alpha$ with $\alpha\not=0$ tends to zero. However the value of the second factor
$q^{L(\alpha)}$ becomes increasingly high when the coordinates of $\alpha$ grow.

As an example, consider the term of (\ref{formula}) corresponding to $\alpha=(r,0,0,r)$. Then $F_\alpha = {\binom n {r,\, r}}^{-1}$, $L(\alpha) = 2\binom r 2$ and\footnote{In this paper the symbol $a_n\sim b_n$ means  that the sequences $a_nb_n^{-1}$ and $a_n^{-1}b_n$ are bounded.}
$$T_\alpha = F_\alpha \cdot q^{L(\alpha)} =  \frac{1}{\binom n {r, \, r}}q^{2\binom r 2}
 \, \sim\,  \left[{\binom n r} p^{\binom r 2}\right]^{-2}.$$
Thus we obtain
\begin{eqnarray}\label{26}
r^2T_{(r,0,0,r)} =o(1),
\end{eqnarray}
by Lemma \ref{lemma4}.

As another example consider the term with $\alpha=(r,0,0,0)$. Then
$$T_\alpha = F_\alpha q^{L(\alpha)} = \frac{\binom {n-2r} r}{\binom n {r, \, r}}q^{\binom r 2}\sim \left[\binom n r p^{\binom r 2}\right]^{-1}.$$
In this case we have \begin{eqnarray}\label{27}
r T_{(r,0,0,0)} =o(1),
\end{eqnarray}
by Lemma \ref{lemma4}.

The term $T_\alpha$ with $\alpha = (1,0,0,0)$ satisfies
\begin{eqnarray}
T_{(1,0,0,0)} \le \frac {r^2}{n}
\end{eqnarray}
as one easily checks.

Next we consider $T_\alpha$ with $\alpha=(r-1, 0,0,0)$.
One has
\begin{eqnarray*}
T_\alpha= r\cdot \frac{\binom{n-2r}{1, \, r}}{\binom n {r,\, r}}q^{\binom {r-1} 2}\sim r(n-3r) \frac{\binom {n-2r} r}{\binom n {r, \, r}}q^{\binom {r-1} 2}\\
\sim \frac{r(n-2r)}{\binom n r} q^{\binom {r-1} 2} \sim np^{r-1}\cdot \left[\frac{r}{{\binom n r}p^{\binom r 2}}\right]\\ \le C'np^{r-1}\, \sim \, C\frac{{\log_q^2 n}}{n}
\end{eqnarray*}
for some constants $C, C'$; here we have used Lemma \ref{lemma4}. Thus, we have the inequality
\begin{eqnarray}
T_{(r-1, 0,0,0)} \le C \cdot \frac{{\log_q^2 n}}{n}.
\end{eqnarray}
Using similar arguments one obtains
\begin{eqnarray}
T_{(r-1, 0,0,r-1)} \le C'' \cdot \frac{{\log_q^4 n}}{n^2},
\end{eqnarray}
where $C''$ is a constant independent of $n$.

As a summary of the above discussion of examples we can make the following claim which will be referred to later:

{\it If $\alpha$ is either $(1,0,0,0)$, or $(r-1, 0,0,0)$, or
$(r-1, 0,0, r-1)$ then }
\begin{eqnarray}\label{r4}
r^4T_\alpha =o(1).
\end{eqnarray}

Recall that $$r= \lfloor 2\log_q n - 2\log_q \log_q n + 2 \log_q(e/2) +1 - \epsilon\rfloor,$$
 and in particular $r\le 2\log_q n$.
Fix $\lambda$ satisfying the inequality
\begin{eqnarray}\label{lambda}
0<\lambda<\frac{1}{1+eq}
\end{eqnarray}
and split the set of all integers in $[0,r]$ into three subsets
$$S_\lambda=\{x\in \N; 0\le x\le (1-\lambda)\log_qn\},$$
$$I_\lambda=\{x\in \N; (1-\lambda)\log_qn < x< (1+\lambda)\log_qn\},$$
$$L_\lambda=\{x\in \N; (1-\lambda)\log_qn \le x\le r\}.$$
Integers lying in $S_\lambda$, $I_\lambda$, and $L_\lambda$ will be called {\it \lq\lq small\rq\rq,\, \lq\lq intermediate\rq\rq\,} and {\it \lq\lq large\rq\rq,\,} correspondingly.

Suppose that $\alpha'\in D$ is obtained from $\alpha=(a, b,c,d)\in D$ by increasing of one of the coordinates by $1$, say, $\alpha'=(a+1,b,c,d)$. Then
the ratio of the corresponding terms of sum (\ref{formula}) equals
$$\frac{T_{\alpha'}}{T_{\alpha}}= \frac{(r-a-b)(r-a-c)}{(a+1)(n-4r+\ell +1)}\cdot q^a,$$
where $\ell =\ell(\alpha) = a+ b+c+ d$. Clearly, one has
$$n/2 \le n-4r+\ell +1\le n,$$ assuming that $n$ is large enough. Hence we obtain
\begin{eqnarray}\label{between}
A\cdot q^a \, \le \, \frac{T_{\alpha'}}{T_{\alpha}}\, \le \, 2\cdot A\cdot q^a
\end{eqnarray}
where
\begin{eqnarray}\label{betweena}
A= \frac{(r-a-b)(r-a-c)}{(a+1)n}.
\end{eqnarray}

If $a\in S_\lambda$ is small then $q^a\le n^{1-\lambda}$, \, $A \le \frac{r^2}{n}$ and $$Aq^a\le \frac{r^2}{n^\lambda}$$ tends to zero as $n\to \infty$.
Hence the ratio which appears in (\ref{between}) is less than $1$ for $n$ large enough.

If $a\in L_\lambda$ is large then $q^a\ge n^{1+\lambda}$, \, $A\ge \frac{1}{2n\log_q n}$ and hence $$Aq^a\ge \frac{n^\lambda}{2\log_q n}$$
tends to infinity for $n\to \infty$. This gives the following statement:
\begin{lemma}\label{lm1}
There exists a constant $N>0$ such that for all $n\ge N$ the following is true:
(1) If $\alpha'\in D$ is obtained from $\alpha\in D$ by adding $1$ to one of its coordinates which is small (see above) then
\begin{eqnarray}
T_\alpha >T_{\alpha'}.
\end{eqnarray}
(2) If $\alpha'\in D$ is obtained from $\alpha\in D$ by adding $1$ to one of its coordinates which is large then
\begin{eqnarray}
T_\alpha <T_{\alpha'}.
\end{eqnarray}
\end{lemma}

\begin{figure}[h]
\begin{center}
\resizebox{7cm}{4cm}
{\includegraphics[21,420][545,767]{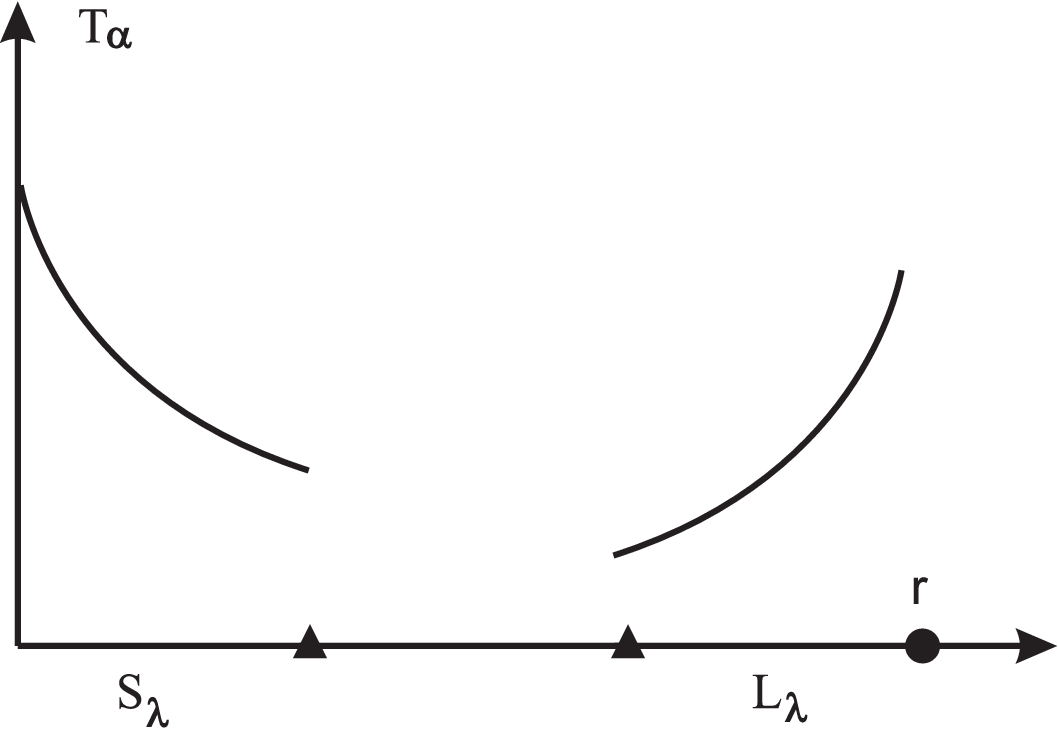}}
\end{center}
\caption{Schematic representation of behavior of $T_\alpha$ with respect to small $a\in S_\lambda$ and large $a\in L_\lambda$ coordinates.} \label{fig2}
\end{figure}

Figure \ref{fig2} illustrates Lemma \ref{lm1}.
Next we analyze the case when one increases an intermediate index.

\begin{lemma}\label{lm2} There exists a constant $N>0$ such that for all $n\ge N$ the following is true:
Suppose that $\alpha'=(a+1,b,c,d)\in D$ is obtained from $\alpha=(a,b,c,d) \in D$ by adding $1$ to one of its coordinates.
 If $a\le r/2$ and ${\mathbf m}(b,c)\not\in S_\lambda$, then
\begin{eqnarray}\label{talpha}
T_\alpha >T_{\alpha'}.
\end{eqnarray}\end{lemma}

\begin{proof} Without loss of generality we may assume that $a\in I_\lambda$ since the case $a\in S_\lambda$ is covered by Lemma \ref{lm1}. Then our assumptions imply that ${\mathbf m}(b,c)\in I_\lambda$, and therefore by symmetry we may assume that $b\in I_\lambda$. Our goal is to estimate the value of $A$ given by (\ref{betweena}). We have
$$a+b>2(1-\lambda) \log_qn$$
and since $r<2\log_qn$ we obtain
\begin{eqnarray}
r-a-b< 2\lambda\log_qn
\end{eqnarray}
and thus the numerator in (\ref{betweena}) satisfies
$$(r-a-b)(r-a-c) < 4 \lambda \log_q^2 n.$$
To estimate the denominator we observe that $a<r/2$ implies
$$q^a \le \frac{eq}{2}\cdot \frac{n}{\log_qn}.$$
Since $a+1 \ge
(1-\lambda) \log_qn$ we obtain
\begin{eqnarray}
2Aq^a \le \frac{4\lambda \log_q^2 n}{(1-\lambda) \log_qn} \cdot \frac{eq}{2}\cdot \frac{n}{\log_qn}\cdot \frac{1}{n}= \frac{4\lambda}{1-\lambda}\cdot \frac{eq}{2}<1;
\end{eqnarray}
the last inequality uses our assumption (\ref{lambda}). This completes the proof of statement (\ref{talpha}). \end{proof}

\begin{lemma}\label{lm3} For $n$ sufficiently large and  $\alpha =(a,0,0,d)\in D$ with $1\le a \le r-1$, one has
\begin{eqnarray}\label{ineqt}
T_\alpha\, \le \, \max\{T_{(1,0,0,d)}, T_{(r-1,0,0,d)}\}.
\end{eqnarray}
\end{lemma}
\begin{proof} The assertion of the Lemma follows from Lemma \ref{lm1} in the case when either $a\in S_\lambda$ or $a\in L_\lambda$. Hence we may assume below that $\alpha=(a, 0,0,d)$ where $a\in I_\lambda$.

Denote
$\alpha'=(a+1, 0,0,d)$ and $\alpha''= (a+2, 0,0,d)$. Then
$$\frac{T_{\alpha''}T_\alpha}{T_\alpha^2} =\left(\frac{r-a-1}{r-a}\right)^2\cdot \left(\frac{a+1}{a+2}\right)\cdot \frac{n-4r+a+d+1}{n-4r+a+d+2}\cdot q.$$
In the RHS of this formula two bracketed factors tend to $1$ as $n\to \infty$; besides $q>1$ . Hence for  $n>N$ large enough one has
\begin{eqnarray}
\frac{T_\alpha T_{\alpha''}}{T_{\alpha'}^2}>1.
\end{eqnarray}
This proves that $\log_q (T_\alpha)$ is convex as function of $a\in I_\lambda$. By Lemma \ref{lm1} this function increases
for $a\in S_\lambda$ and decreases for $a\in L_\lambda$. This implies (\ref{ineqt}).
\end{proof}

Now we are able to complete the proof of Theorem \ref{main1}. Recall that we have to show that the sum $\sum_{\alpha\in D'}T_\alpha$ tends to $0$ as $n\to \infty$ where $D'=D-\{(0,0,0,0)\}$. Consider the subset $\tilde D\subset D$ consisting of vectors with at least one coordinate equal $r$. Each $\alpha\in \tilde D$ has the form $\alpha=(r,0,0,d)$ (up to symmetry) where $d=0, \dots, r$. Applying Lemma \ref{lm3} we obtain that
$$T_\alpha \le \max\{T_{(r, 0,0,0)}, T_{(r,0,0,r)}\}.$$
Since the cardinality of $\tilde D$ does not exceed $5r$, we obtain, using (\ref{26}) and (\ref{27}), that
\begin{eqnarray}\label{41}
\sum_{\alpha\in \tilde D} T_\alpha = o(1).
\end{eqnarray}

Each vector $\alpha\in D'$ may have at most two large coordinates. Decompose $$D'- \tilde D=D'_0\cup D'_1\cup D'_2,$$
where $D'_i$ denotes the set all vectors in $\tilde D$ having exactly $i$ large coordinates, $i=0,1,2$.

Suppose that $\alpha \in D_2'$. Without loss of generality we may assume that $a$ and $d$ are large and $b$ and $c$ are small, i.e. $a,d\in L_\lambda$, $ b,c\in S_\lambda$.
Applying Lemma \ref{lm1} we obtain
$T_\alpha \leq T_{(a,0,0,d)}.$
Since $a\not= r\not= d$ we may engage Lemma \ref{lm3} to obtain
\begin{eqnarray}\label{42}
T_\alpha \le \max\{T_{(1,0,0,r-1)}, T_{(r-1,0,0,r-1)}\}.
\end{eqnarray}
Now, taking into account (\ref{26}), (\ref{27}) and (\ref{r4}), we obtain
\begin{eqnarray}\label{two}
\sum_{\alpha\in D'_2} T_\alpha =o(1).
\end{eqnarray}

Consider now the sum $\sum_{\alpha\in D'_1}T_\alpha$. In this case the vector $\alpha=(a, b,c,d)$ contains one large index. Assume that $a$ is large. Then $b, c$ must be small and applying Lemma \ref{lm1} and Lemma \ref{lm3} we obtain
\begin{eqnarray*}
\qquad T_\alpha\le T_{(a,0,0,d)} \le T_{(r-1,0,0,d)}\le \max\{T_{(r-1,0,0,0)}, T_{(r-1,0,0,r-1)}\}.
\end{eqnarray*}
Now (\ref{r4}) implies that
\begin{eqnarray}\label{one}
\sum_{\alpha\in D_1'}T_\alpha =o(1).
\end{eqnarray}

Next we show that for any $\alpha\in D'_0$ one has
\begin{eqnarray}\label{four4}
T_\alpha \le \max\{T_{(1,0,0,0)}, T_{(r-1,0,0,0)}, T_{(r-1,0,0,r-1)}\}
\end{eqnarray}
which in view of (\ref{r4}) would imply that
\begin{eqnarray}\label{zero}
\sum_{\alpha\in D'_0}T_\alpha =o(1).
\end{eqnarray}
The combination of (\ref{41}), (\ref{two}), (\ref{one}) and (\ref{zero}) give Theorem \ref{main1}.

To prove (\ref{four4}) consider $\alpha=(a,b,c,d)\in D'_0$. Note that coordinates $a, b, c, d$ can be either small or intermediate. Assume first that all coordinates
$a, b, c,d$ are small. Then $T_\alpha\le T_{(1,0,0,0)}$ (by Lemma \ref{lm1}) implying (\ref{four4}).

Suppose now that exactly one of the coordinates of $\alpha$ is intermediate. If $a$ is intermediate and $b, c,d$ are small then
$$T_\alpha\le T_{(a,0,0,0)}\le \max\{T_{(1,0,0,0)}, T_{(r-1,0,0,0)}\}$$
(by Lemma \ref{lm1} and Lemma \ref{lm3}) proving (\ref{four4}).

Suppose now that two coordinates of $\alpha$ are intermediate. Taking into account symmetry (the action of $G$ on $D$, see above), this case can be subdivided into two subcases:
(i) $a$ and $b$ are intermediate and (ii) $a$ and $d$ are intermediate. In the subcase (i), since $a+b\le r$, either $a\le r/2$, or $b\le r/2$ and we may apply Lemma \ref{lm2}.
Assuming that $a\le r/2$ we obtain
$$T_\alpha \le T_{(0,b,0,0)} \le \max\{T_{(1,0,0,0)}, T_{(r-1,0,0,0)}\},$$
implying (\ref{four4}). In the subcase (ii), we know that $b, c$ are small hence $T_\alpha\le T_{(a,0,0,d)}$ and application of Lemma \ref{lm3} gives (\ref{four4}).

In the remaining case when $\alpha\in D'_0$ has three or four intermediate indices we know that at least two of these indices are $\le r/2$ and by Lemma \ref{lm2}
one has $$T_\alpha\le T_{\alpha'}$$ where $\alpha'$ is obtained from $\alpha$ by replacing by zeros two coordinates which were $\le r/2$. To estimate $T_{\alpha'}$ one applies Lemma \ref{lm3} leading again to (\ref{four4}).
This completes the proof of Theorem \ref{main1}. \qed


\begin{thebibliography}{99}


  \bibitem {AS} N. Alon, J. Spencer, \textit{The probabilistic method}. John Wiley \& Sons, Inc., Hoboken, NJ, 2008.


\bibitem{BHK} E. Babson, C. Hoffman and M. Kahle, \textit{The fundamental group of random $2$-complexes}, 
J. Amer. Math. Soc. {\bf 24} (2011), 1-28. 


  \bibitem{BE} B. Bollob\'{a}s and P. Erd\H{o}s, \textit{Cliques in random graphs}, Mathem. Proceedings of Cambridge Phil. Soc. {\bf 80}(1976), 419 - 427.

  \bibitem{B} B. Bollob\'{a}s, \textit{Random Graphs}, Cambridge University Press, 2008.

\bibitem{B1} B. Bollob\'{a}s, \textit{Random graphs}. In {\it Combinatorics}, Proceedings, Swansea 1981, London Math. Soc. Lecture Note Ser. 52, Cambridge University Press, Cambridge, 80-102. 


  \bibitem{C} R. Charney,
\textit{An introduction to right-angled Artin groups},
Geom. Dedicata 125 (2007), 141--158.

\bibitem{ChF} R. Charney and M. Farber, \textit{Random groups arising as graph products}, arXiv:1006.3378.

\bibitem{CFK} A. Costa, M. Farber and T. Kappeler, \textit{Topology of random 2-complexes}, arXiv:1006.4229.



\bibitem{CF} A. Costa, M. Farber, \textit{Motion planning in spaces with small fundamental groups}, Communications in Contemporary Mathematics, to appear.

\bibitem{CP} D. Cohen and G. Pruidze,  \textit{Motion planning in tori}, Bull. Lond. Math. Soc. 40 (2008), no. 2, 249--262.

\bibitem{CW} J. Crisp and B. Wiest, \textit{Embeddings of graph braid and surface groups in right-angled Artin groups and braid groups.}
Algebr. Geom. Topol. 4 (2004), 439--472. 

\bibitem{ER} P.\ Erd\H{o}s, A.\ R\'enyi, {On the evolution of 
random graphs}, Publ.\ Math.\ Inst.\ Hungar.\ Acad.\ Sci.\ {\bf 5}
(1960), 17--61.

\bibitem{F1} M. Farber,
\textit{Topological complexity of motion planning,} Discrete and Comput. Geom., \textbf{29}(2003), 211--221.

\bibitem{F2} M. Farber,
\textit{Instabilities of Robot Motion}, Topology and its applications, \textbf{140}(2004), 245-266.

\bibitem{FTY} M. Farber, S. Tabachnikov and S. Yuzvinsky,
\textit{Topological robotics: motion planning in
projective spaces}, Int. Math. Res. Not., \textbf{34}(2003), 1853--1870.


\bibitem{invit} M. Farber, \textit{Invitation to topological robotics}, Zurich Lectures in Advanced Mathematics, EMS, 2008.


\bibitem{F1} M. Farber, 
\emph{Topology of random linkages, }
Algebraic and Geometric Topology, 8(2008), 155 - 171.


\bibitem {FK} M. Farber and T. Kappeler,
\textit{Betti numbers of random manifolds}, 
Homology, Homotopy and Applications, Vol. 10 (2008), No. 1, pp. 205 - 222.

\bibitem{Gr} M. Gromov, \emph{Asymptotic invariants of infinite groups}, LMS, 1993.

\bibitem{G} M. Gromov, \textit{Random walk in random groups},  Geom. Funct. Anal. 13, No. 1, 73-146 (2003).


  \bibitem{JLR} S. Janson, T. Luczak, A. Rucinski, \textit{Random graphs}, John Wiley and Sons, 2000.

\bibitem{KS} I. Kapovich, P. Schupp, \textit{On group-theoretic models of randomness and genericity}, Groups, Geom. Dyn. 2(2008), No. 3, 383-404.

\bibitem{KR} M. Karo\'{n}ski, A. Ruci\'{n}ski, \textit{On the number of strictly balanced subgraphs of a random graph}, In {\it Graph Theory}, Proceedings, 
Lag\'ow, 1981, M. Borowiecki et al editors, Lecture Notes in Math. 1018(1983), 79-83. 

  \bibitem{La} J.-C. Latombe, \textit{Robot Motion Planning}, Kluwer, Dordrecht, 1991.

\bibitem{LM} N.\ Linial, R.\ Meshulam, {Homological connectivity
  of random $2$-complexes}, Combinatorica {\bf 26} (2006),  475--487.


  \bibitem{M1} D.W. Matula, \textit{On the complete subgraphs in a random graph}, Combinatory Mathematics and its Applications, Chapel Hill, North Carolina, 1970, pp. 356 - 369.

  \bibitem{M2} D.W. Matula, \textit{The largest clique size in a random graph}, Tech. Rep., Dept. Comput. Sci., Southern Methodist University, Dallas, 1976.

  \bibitem{MV} J. Meier, L. VanWyk,
\textit{The Bieri-Neumann-Strebel invariants for graph groups},
Proc. London Math. Soc. (3) 71 (1995), no. 2, 263--280.



\bibitem{MW} R.\ Meshulam, N.\ Wallach, {Homological
  connectivity of random $k$-complexes}, Random Structures \& Algorithms 
  \textbf{34} (2009), 408--417. 

\bibitem{Oll} Y. Ollivier, \textit{A January 2005 invitation to random groups}. Ensaios Matemáticos, 10. Sociedade Brasileira de Matemática, Rio de Janeiro, 2005. ii+100 pp. 

\bibitem{Sch} K. Sch\"urger, \textit{Limit theorems for complete subgraphs of random graphs}, Per. Math. Hungar, {\bf{10}}, 47 - 53.


\bibitem{S} A. Shalev, \textit{Probabilistic group theory and Fuchsian groups}. Infinite groups: geometric, combinatorial and dynamical aspects, 363--388, 
Progr. Math., 248, Birkhäuser, Basel, 2005. 


\bibitem{Zuk} A. Zuk, \textit{Property (T) and Kazhdan constants for discrete groups.} 
Geom. Funct. Anal. {\bf 13} (2003), no. 3, 643--670. 


  \end{thebibliography}
\end{document}